\documentclass[twoside,a4paper,12pt,centertags,reqno]{amsart} % 13pt displays better than 12pt
\usepackage{amsmath,amssymb,verbatim,vmargin}
\usepackage{color}
\usepackage{tikz}
\usepackage{hyperref}
\hypersetup{colorlinks,bookmarks=true,linktocpage=true,citecolor=blue,linkcolor=magenta}

\usepackage{mathptmx}   % based on the Times Roman font

\allowdisplaybreaks % align break
\usepackage{color}
\usepackage[all]{xy} % diagrams

\usepackage{mathtools}
\usepackage{MnSymbol} %wedge righthalfcup

\DeclareFontFamily{U}{mathx}{\hyphenchar\font45}
\DeclareFontShape{U}{mathx}{m}{n}{
      <5> <6> <7> <8> <9> <10>
      <10.95> <12> <14.4> <17.28> <20.74> <24.88>
      mathx10
      }{}
\DeclareSymbolFont{mathx}{U}{mathx}{m}{n}
\DeclareFontSubstitution{U}{mathx}{m}{n}
\DeclareMathAccent{\widecheck}{0}{mathx}{"71}
\DeclareMathAccent{\wideparen}{0}{mathx}{"75}

\theoremstyle{plain}
\newtheorem{thm}{Theorem}[section]

\theoremstyle{definition}

\newtheorem{rem}[thm]{Remark}

\newtheorem{obs}[thm]{Observation}

\usepackage{enumerate}

\newcommand{\bR}{{\mathbb R}}

\newcommand{\bN}{{\mathbb N}}

\newcommand{\cA}{{\mathcal A}}

\newcommand{\cD}{{\mathcal D}}

\def\barint_#1{\mathchoice
            {\mathop{\vrule width 6pt
height 3 pt depth -2.5pt
                    \kern -9.5pt
\intop \kern -4pt}\nolimits_{#1}}%
            {\mathop{\vrule width 5pt height
3 pt depth -2.6pt
                    \kern -6.5pt
\intop \kern -4pt}\nolimits_{#1}}%
            {\mathop{\vrule width 5pt height
3 pt depth -2.6pt
                    \kern -6pt
\intop \kern -4pt}\nolimits_{#1}}%
            {\mathop{\vrule width 5pt height
3 pt depth -2.6pt
          \kern -6pt \intop \kern -4pt}\nolimits_{#1}}}
          
           \def\bariint_#1{\mathchoice
            {\mathop{\vrule width 15pt
height 3 pt depth -2.5pt
                    \kern -15.8pt
\intop \kern -8pt\intop \kern -4pt}\nolimits_{#1}}%
            {\mathop{\vrule width 9pt height
3 pt depth -2.6pt
                    \kern -10.5pt
\intop \kern -8pt\intop \kern -4pt}\nolimits_{#1}}%
            {\mathop{\vrule width 9pt height
3 pt depth -2.6pt
                    \kern -10pt
\intop \kern -8pt\intop \kern -4pt}\nolimits_{#1}}%
            {\mathop{\vrule width 9pt height
3 pt depth -2.6pt
          \kern -8pt \intop \kern -10pt\intop \kern -4pt}
      \nolimits_{  #1}}}

\def\barintlim_#1{\mathchoice
            {\mathop{\vrule width 6pt
height 3 pt depth -2.5pt
                    \kern -8.8pt
\intop \kern -4pt}\limits_{#1}}%
            {\mathop{\vrule width 5pt height
3 pt depth -2.6pt
                    \kern -6.5pt
\intop \kern -4pt}\limits_{#1}}%
            {\mathop{\vrule width 5pt height
3 pt depth -2.6pt
                    \kern -6pt
\intop \kern -4pt}\limits_{#1}}%
            {\mathop{\vrule width 5pt height
3 pt depth -2.6pt
          \kern -6pt \intop \kern -4pt}\limits_{#1}}}
          
           \def\bariintlim_#1{\mathchoice
            {\mathop{\vrule width 15pt
height 3 pt depth -2.5pt
                    \kern -15.8pt
\intop \kern -8pt\intop \kern -4pt}\limits_{#1}}%
            {\mathop{\vrule width 9pt height
3 pt depth -2.6pt
                    \kern -10.5pt
\intop \kern -8pt\intop \kern -4pt}\limits_{#1}}%
            {\mathop{\vrule width 9pt height
3 pt depth -2.6pt
                    \kern -10pt
\intop \kern -8pt\intop \kern -4pt}\limits_{#1}}%
            {\mathop{\vrule width 9pt height
3 pt depth -2.6pt
          \kern -8pt \intop \kern -10pt\intop \kern -4pt}
      \limits_{  #1}}}
          
\renewcommand{\iint}{\int \kern -3pt\int}       

\numberwithin{equation}{section}
\makeatletter
\@namedef{subjclassname@2020}{\textup{2020} Mathematics Subject Classification}
\makeatother

\newcommand{\nocontentsline}[3]{}
\let\origcontentsline\addcontentsline
\newcommand\stoptoc{\let\addcontentsline\nocontentsline}
\newcommand\resumetoc{\let\addcontentsline\origcontentsline}

\title[Multiple mass drop for the Davey--Stewartson II equation]{On the Moutard transformation singularity for the Davey--Stewartson II equation}

\author{Yi C. Huang} 
\address{School of Mathematical Sciences, Nanjing Normal University, Nanjing 210023, People's Republic of China}
\email{Yi.Huang.Analysis@gmail.com}
\urladdr{https://orcid.org/0000-0002-1297-7674}

\date{\today} 

\thanks{Research of the author is partially supported by the National NSF grant of China (no. 11801274), the JSPS Invitational Fellowship for Research in Japan (no. S24040),
and the Open Projects from Yunnan Normal University (no. YNNUMA2403) and Soochow University (no. SDGC2418).
This text is written during a visit to Novosibirsk and the author thanks sincerely Iskander A. Ta{\u\i}manov for enlightening discussions.}

\begin{document}

\begin{abstract}
We indicate explicitly how to obtain the Moutard transformation singularity for the Davey--Stewartson II equation with two indeterminancy points.
We also justify an arbitrary order ``mass drop" phenomenon via this singularity formation scheme.
\end{abstract}

\maketitle

%\tableofcontents

\section{Introduction: the Davey--Stewartson II equation}

In this Letter, following Ta{\u\i}manov \cite{Tai21, Tai24}, 
we shall consider the Moutard transformation solutions of the following Davey--Stewartson II equation \cite{DavSte74}:
\begin{equation}
\label{ds2}
U_t = i\left(U_{zz}+U_{\overline{z}\,\overline{z}} + 2(V+\overline{V})U\right), \ \ \ V_{\overline{z}} = (|U|^2)_z, % position of the factor 2 is harmless; in accordance with thm, better in first equation
\end{equation}
which is the compatibility condition for the linear problems
$$
\cD \Psi =0 , \ \ \ \partial_t \Psi = \cA\Psi,
$$
where
$$
\cD= \left(
\begin{array}{cc}
0 & \partial \\
-\overline{\partial} & 0
\end{array}
\right) + \left(
\begin{array}{cc}
U & 0 \\
0 & \overline{U}
\end{array}
\right)
$$
is a two-dimensional Dirac operator with a complex-valued potential $U$ with
$$\partial = \frac{1}{2}\left(\frac{\partial}{\partial x} - i\frac{\partial}{\partial y}\right),
\quad\overline{\partial} = \frac{1}{2}\left(\frac{\partial}{\partial x} + i\frac{\partial}{\partial y}\right),$$
the time-dependent generator $\cA$ is given by
$$\cA = i\left(\begin{array}{cc} -\partial^2 - V & \overline{U}\,\overline{\partial} - \overline{U}_{\overline{z}} \\
U\partial - U_z & \overline{\partial}^2 + \overline{V} \end{array}\right), $$
and the solution $\Psi$ is taken as the following matrix form 
$$
\Psi =
\left(\begin{array}{cc} \psi_1 & -\overline{\psi}_2 \\
\psi_2 & \overline{\psi}_1 \end{array} \right).
$$

Note that if $(U_0,V_0)$ solves \eqref{ds2}, then $(\overline U_0,V_0)$ solves
$$U_t = -i\left(U_{zz}+U_{\overline{z}\,\overline{z}} + 2(V+\overline{V})U\right), \ \ \ V_{\overline{z}} = (|U|^2)_z.$$
Hence, \eqref{ds2} is also the compatibility condition of the ``conjugate" problems
$$\cD^\vee \Phi = 0, \ \ \Phi_t = \cA^\vee \Phi,$$
where
$$\cD^\vee = \left(\begin{array}{cc} 0 & \partial \\
-\overline{\partial} & 0 \end{array}\right) +
\left(\begin{array}{cc} \overline{U} & 0 \\
0 & U \end{array}\right),$$
and
$$\cA^\vee= - i\left(\begin{array}{cc} -\partial^2 - V & U\overline{\partial} - U_{\overline{z}} \\
\overline{U}\partial - \overline{U}_z & \overline{\partial}^2 + \overline{V} \end{array}\right).$$
Analogously, the solution $\Phi$ is taken as the following matrix form 
$$
\Phi = \left(\begin{array}{cc} \varphi_1 & -\overline{\varphi}_2 \\
\varphi_2 & \overline{\varphi}_1 \end{array} \right).
$$

If $U=0$, $(\psi_1,\varphi_1)$ is pair of holomorphic functions, while $(\psi_2,\varphi_2)$ anti-holomorphic.
When $U=\overline{U}$, the ``conjugate" problems are irrelevant and one regards $\Psi=\Phi$.
The pair $(\Psi,\Phi)$ serves as a Davey--Stewartson II deformation of surfaces in $\bR^4$, see e.g. \cite{Tai22}.

\section{Moutard singularity on two points or multiple zeroes}

Using the Moutard transformation of the solution $U=V=0$ with spinors
$$
\Psi = \left(\begin{array}{cc} g' & 0 \\ 0 & \overline{g'} \end{array}\right), \ \
\Phi = \left(\begin{array}{cc} f^\prime & i \\ i & \overline{{f}^\prime} \end{array}\right),
$$
the following two-parameter family of solutions were found in \cite[Theorem 4.1]{Tai24}.

\begin{thm}[Ta{\u\i}manov, 2024] \label{twopara}
Let $f(z,t)$ and $g(z,t)$ be two functions which are holomorphic in $z$ and satisfy the equations
\begin{equation} \label{fg}
\frac{\partial f}{\partial t} = i\frac{\partial^2 f}{\partial z^2},\quad \frac{\partial g}{\partial t} = -i\frac{\partial^2 g}{\partial z^2}.
\end{equation}
Let $h$ be a function which is holomorphic in $z$ and satisfies the relations
$$h'=f'g',\quad ih_t=g''f'-g'f''.$$ 
Then
\begin{equation} \label{Tai}
U = i\overline{g'}\frac{f^\prime g - h}{|g|^2 + | h|^2}
\end{equation}
satisfies the Davey--Stewartson II equation \eqref{ds2}.
\end{thm}

It is noteworthy to point out that in the singularity formation scheme of Moutard transformation,
the singularity points $(z_0,t_0)$ are encoded in the zero set of $|g|^2 + |h|^2$. 

\begin{rem}
Let us first write down the example given in \cite{Tai24}.
Take $$f(z,t)=z^2+2it,\quad g(z,t)=z^3-6itz,$$
and choose
$$h=\frac32 z^4-6itz^2+6t^2.$$
We derive the following solution
$$\begin{aligned}
U &=\frac{i\left(3\overline z^2+6it\right)\left(\frac12 z^4-6itz^2-6t^2\right)}{|z^3-6itz|^2 + \left| \frac32 z^4-6itz^2+6t^2\right|^2}\\
&=O(|z|^{-2}),\qquad |z|\rightarrow+\infty.
\end{aligned}$$
For $t=0$ it has the indeterminacy type singularity $\sim\frac32i\frac{z^2}{|z|^2}$ as $|z|\rightarrow0$. 
\end{rem}

The first observation of this Letter can be summarized as follows:

\begin{obs} \label{2indeter}
In above example, we modify the data set $(f,g)$ in the following way
$$f(z,t)=(z^2+2it)-z,\quad g(z,t)=(z^3-6itz)-z.$$
The natural observation (or rather, motivation) here is that the lower order modification terms is set to solve the \textit{linear} equations for $f$ and $g$.
Thus it remains to compute
$$f'(z,t)=2z-1,\quad g'(z,t)=3z^2-(1+6it)$$
and
$$h'(z,t)=f'(z,t)g'(z,t)=6z^3-3z^2-2(1+6it)z+(1+6it).$$
We choose
$$h(z,t)=\frac32 z^4-z^3-(1+6it)z^2+(1+6it)z-i(t^2+4it^3)-\frac12.$$
Here, by adding the temporal term $-i(t^2+4it^3)$ we meet the constraint $ih_t=g''f'-g'f''$, 
and by adding the constant $-\frac12$, we see that $h|_{t=0}$ vanishes at $z=\pm1$.
By Theorem \ref{twopara}, the data triple $(f,g,h)$ yields the following solution for the Davey--Stewartson II equation
$$\begin{aligned}
U &=i\overline{g'}\frac{f^\prime g - h}{|g|^2 + | h|^2}\\
&= \frac{i\left(3\overline z^2-(1-6it)\right)\left(\frac12 z^4-(1+6it) z^2+i(t^2+4it^3)+\frac12\right)}{|(z^3-6itz)-z|^2 + \left| \frac32 z^4-z^3-(1+6it)z^2+(1+6it)z-i(t^2+4it^3)-\frac12\right|^2}.
\end{aligned}$$
Note that
$$h|_{t=0}=\frac32 z^4-z^3-z^2+z-\frac12=(z^2-1)\left(\frac32 z^2-z+\frac12\right),$$
$$g|_{t=0}=z(z^2-1),\quad g'|_{z=\pm1}\neq0,$$
and
$$(f'g-h)|_{t=0}=\frac{(z^2-1)^2}{2}.$$
Thus the solution $U$ at $t=0$ has a two-point indeterminancy type singularity 
$$U|_{t=0}\sim C_{\pm1}\frac{(z^2-1)^2}{|z^2-1|^2},\quad z\rightarrow \pm1,$$
The involved implicit constants $C_{\pm1}$ can be computed directly, and we omit this step.
\end{obs}

\begin{rem}
For higher order polynomials one can construct in the similar manner the multiple-point indeterminancy singularity.
We choose not to dwell on that issue.
\end{rem}

Note that the solution in Observation \ref{2indeter} has mass $\|U(\cdot,t)\|_2^2\equiv 4\pi$,
except at the singular time $t=0$, it has a ``mass drop" with magnitude $2\pi$ (with $2$ being the number of indeterminancy points).
The second observation of this Letter can be summarized as follows:

\begin{obs}
According to \eqref{fg}, we construct for $2\leq m\in\bN$ and $1\leq n\in\bN$ that
$$f_m=z^m+im(m-1)tz^{m-2}+\cdots,$$ 
$$g_n=z^n-in(n-1)tz^{n-2}-\cdots,$$
and for convenience, with the normalization $f_m(0,0)=g_n(0,0)=0$. 
For example, we have
$$\begin{aligned}
f_2(z,t)&=z^2+2it, \\
f_3(z,t)&=z^3+6itz,\\
f_4(z,t)&=z^4+12itz^2-12t^2,\\ 
f_5(z,t)&=z^5+20itz^3-60t^2z
\end{aligned}$$
and
$$\begin{aligned}
g_1(z,t)&=z, \\
g_2(z,t)&=z^2-2it, \\
g_3(z,t)&=z^3-6itz,\\
g_4(z,t)&=z^4-12itz^2-12t^2,\\ 
g_5(z,t)&=z^5-20itz^3-60t^2z.
\end{aligned}$$
In constructing the function $h$ we also assume $h(0,0)=0$,
thus the solution $U$ in \eqref{Tai} with data $(f_m,g_n)$ has mass $\|U(\cdot,t)\|_2^2\equiv (m+n-1)\pi$,
except at the singular time $t=0$, it has a ``mass drop" with magnitude $n\pi$ (with $n$ being the order of $g$) to $(m-1)\pi$.
The normalization on $h$ can not be fulfilled when $m=1$ unless $U$ is the trivial zero solution.
\end{obs}

\bigskip

\stoptoc

\section*{\textbf{Compliance with ethical standards}}

\bigskip

\textbf{Conflict of interest} The author has no known competing financial interests
or personal relationships that could have appeared to influence this reported work.

\bigskip

\textbf{Availability of data and material} Not applicable.

\bigskip

\resumetoc

\bibliographystyle{alpha}

\bibliography{Hua-Moutard-DavSteII-Multiple} 
 
\end{document}